\documentclass[11pt,letterpaper]{amsart}

\usepackage[T1]{fontenc}
\usepackage[utf8]{inputenc}
\usepackage{lmodern}

\usepackage{amsmath,amssymb,amsfonts,amsthm}
\usepackage{geometry}
\usepackage{enumitem}
\usepackage{booktabs}
\usepackage{subcaption}
\usepackage{caption}
\usepackage{xcolor}
\definecolor{wiseblue}{rgb}{0.01,0.12,0.42}
\definecolor{softblue}{RGB}{170,205,230}
\definecolor{sagegreen}{RGB}{177,211,179}
\definecolor{warmcoral}{RGB}{238,154,138}
\definecolor{duckyellow}{RGB}{255,239,153}
\definecolor{circleedge}{RGB}{18,63,112}
\definecolor{contactedge}{RGB}{18,63,112}

\usepackage[hidelinks]{hyperref}

\usepackage{tikz}

\numberwithin{equation}{section}

\newtheorem{theorem}{Theorem}[section]
\newtheorem{proposition}[theorem]{Proposition}
\newtheorem{lemma}[theorem]{Lemma}
\newtheorem{corollary}[theorem]{Corollary}

\theoremstyle{definition}
\newtheorem{definition}[theorem]{Definition}

\theoremstyle{remark}
\newtheorem{remark}[theorem]{Remark}

\DeclareMathOperator{\rank}{rank}
\DeclareMathOperator{\Hess}{Hess}

\newcommand{\R}{\mathbb R}
\newcommand{\Q}{\mathbb Q}
\newcommand{\D}{\mathcal D}

\newcommand{\Roll}{\mathrm{\mathbf{Roll}}}
\newcommand{\Adm}{\mathrm{\mathbf{Adm}}}

\newcommand{\bc}{\mathbf c}
\newcommand{\bv}{\mathbf v}
\newcommand{\bw}{\mathbf w}
\newcommand{\bx}{\mathbf x}
\newcommand{\bs}{\mathbf s}
\newcommand{\bz}{\mathbf z}
\newcommand{\bu}{\mathbf u}

\newcommand{\blambda}{\boldsymbol\lambda}

\newcommand{\cbar}{\overline{\mathbf c}}
\newcommand{\vbar}{\overline{\mathbf v}}
\newcommand{\wbar}{\overline{\mathbf w}}

\newcommand{\Euc}{\mathrm{E}(2)}

\begin{document}

\title{Second Moment Penny Packings via Rolling Spaces}

\author{Jos\'e Ayala Hoffmann}
\thanks{ORCID:
\href{https://orcid.org/0000-0003-2649-6676}
{0000-0003-2649-6676}.
}
\address{Universidad de Tarapac\'a, Iquique, Chile}
\email{jayalhoff@gmail.com}

\author{Fabi\'an Henry Vilaxa}
\address{Universidad de Tarapac\'a, Iquique, Chile}
\email{fabian.henry.vilaxa@alumnos.uta.cl}

\subjclass[2020]{Primary 52C15; Secondary 52C25, 68V05}

\keywords{
disk packing,
penny graph,
second moment,
rolling space,
rigidity,
computer-assisted proof,
interval arithmetic
}

\begin{abstract} \baselineskip=20 true pt
\baselineskip=1.1\normalbaselineskip

We solve the unrestricted minimum second moment problem for six unit disks
in the plane by proving that the minimum is \(29/3\) attained uniquely
up to isometry and relabelling by one disk together with five of its
six neighbours in the triangular lattice.

The realised contacts themselves define the space of admissible motions,
while the directions preserving those contacts to first order define
the rolling space consisting of motions whose relative velocities
have zero normal component at every realised contact. 
Together, these spaces lead to an exact tangent cone formula and to a
criterion by which a strictly positive centroidal contact balance,
combined with rolling rigidity, implies strict local minimality
modulo Euclidean motions. Our method recovers the classical five disk minimum \(34/5\) of Graham and Sloane, previously the largest unrestricted case known exactly.

Geometric reductions and a graph census reduce the global search
to finitely many contact graph classes.
Exact analysis and validated interval arithmetic then exclude
competing configurations, establishing the global bounds and uniqueness. Our method allows us to pass from putative lattice extremisers to
unrestricted optimality. 

\end{abstract}

\maketitle

\baselineskip=20 true pt
\baselineskip=1.1\normalbaselineskip
\section{Prelude}
\label{sec:introduction}

How tightly can a finite collection of pennies be packed around its centroid without overlapping? The second moment of the disk centres
measures their spread while penalising outliers. Graham and Sloane studied this problem for the planar case~\cite{grahamsloane}, motivated in part by two dimensional codes of minimum energy. They determined the unrestricted optimum for five
disks and conjectured that, except when \(n=4\), every optimal packing
lies in the triangular lattice, up to isometry.
Chow later on provided an algorithm for finding least second moment packings
within the triangular lattice~\cite{chow}. This result identifies candidates, but the passage from a lattice optimum to an unrestricted optimum is substantial. Optimality must be proved against every feasible planar configuration including configurations outside the triangular
lattice. 

Let \(\bc=(\bc_1,\ldots,\bc_n)\in(\R^2)^n\) be the centres of \(n\)
non overlapping unit disks, and let
\[
 \cbar=\frac1n\sum_{i=1}^n\bc_i
\]
be their centroid. We consider the centroidal second moment
\begin{equation}
\label{eq:objective}
 M(\bc)
 =
 \frac12\sum_{i=1}^n\|\bc_i-\cbar\|^2.
\end{equation}
Equivalently,
\begin{equation}
\label{eq:pairwise-objective}
 M(\bc)
 =
 \frac{1}{2n}
 \sum_{1\le i<j\le n}
 \|\bc_i-\bc_j\|^2.
\end{equation}

Our normalisation differs from that of Graham and
Sloane by a factor of two, $M=2U$,
where $U$ denotes their second moment. 

We solve the unrestricted minimum second moment problem for six unit disks in the plane, proving that the minimum is \(29/3\), attained uniquely up to isometry and relabelling by one disk together with five of its six neighbours in the triangular lattice. Our method also recovers the classical five disk minimum \(34/5\) of Graham and Sloane.

The calculus for realised hard disk configurations originated in our
study of perimeter minimisation for hard disk clusters~
\cite{ayalahoffmannvilaxa-intrinsic}. Here we apply its intrinsic
contact geometry to the centroidal second moment, developing the
corresponding balance and local minimality criteria that underpin
our unrestricted six disk classification. For a hard disk configuration \(\bc=(\bc_1,\ldots,\bc_n)\), the local analysis is organised around its realised contacts. A contact may open but cannot shorten to first order.
This allows us to define the admissible cone \(\Adm(\bc)\).
In fact, this cone coincides exactly with the tangent cone
to the full hard disk configuration space at \(\bc\).
Equality at every contact defines the rolling space
\(\Roll(\bc)\), consisting of the allowed motions that
preserve the contacts to first order. In the perimeter problem, the local analysis must accommodate changes in the realised contacts and in the perimeter functional, which is piecewise smooth (its local formula changes with the combinatorics of the convex hull even within the same contact class). On the other hand, the centroidal second moment \(M\) is smooth and quadratic on the whole ambient configuration space, so the changing local geometry comes exclusively from the realised non overlap constraints. This makes the second moment problem a particularly transparent setting in which to separate the intrinsic feasible geometry from the chosen functional.

This exact tangent cone identity removes the need for a constraint qualification at the first order level. At an ambient local minimiser, \(DM(\bc)\) is nonnegative on \(\Adm(\bc)\), and a polar cone argument yields a centroidal contact balance with nonnegative multipliers. A strictly positive balance together with rolling rigidity implies strict local minimality modulo orientation preserving Euclidean motions without assuming independence of the active contact gradients. Both target configurations, shown in
Figure~\ref{fig:all-candidates}, nevertheless satisfy
the linear independence constraint qualification,
as verified by the contact rank computations.

\begin{figure}[!htbp]
\centering
\renewcommand\thesubfigure{\Alph{subfigure}}

\def\cvone{softblue}               
\def\cvtwo{sagegreen}              
\def\cvthree{warmcoral}            
\def\cvfour{duckyellow}            
\def\cvfive{softblue!55!duckyellow}
\def\cvsix{warmcoral!55!duckyellow} 
\def\cvseven{wiseblue!28}          

\begin{subfigure}[t]{0.48\textwidth}
\centering
\begin{tikzpicture}[scale=1]

\def\rad{1}

\path[use as bounding box] (-3.15,-2.35) rectangle (3.15,2.35);

\coordinate (A) at ({-2},{-2*sqrt(3)/5});
\coordinate (B) at ({0},{-2*sqrt(3)/5});
\coordinate (C) at ({2},{-2*sqrt(3)/5});
\coordinate (D) at ({-1},{3*sqrt(3)/5});
\coordinate (E) at ({1},{3*sqrt(3)/5});

\fill[\cvone]  (A) circle (\rad); 
\fill[\cvtwo]  (B) circle (\rad); 
\fill[\cvthree](C) circle (\rad); 
\fill[\cvfour] (D) circle (\rad); 
\fill[\cvfive] (E) circle (\rad); 

\foreach \P in {A,B,C,D,E}{
    \draw[circleedge,line width=0.7pt] (\P) circle (\rad);
}

\draw[contactedge,line width=0.8pt] (A)--(B)--(C);
\draw[contactedge,line width=0.8pt] (D)--(E);
\draw[contactedge,line width=0.8pt] (A)--(D);
\draw[contactedge,line width=0.8pt] (D)--(B);
\draw[contactedge,line width=0.8pt] (B)--(E);
\draw[contactedge,line width=0.8pt] (E)--(C);

\node[xshift=-5pt,yshift=-5pt] at (A) {\small 1};
\node[yshift=-6pt]             at (B) {\small 2};
\node[xshift=5pt,yshift=-5pt]  at (C) {\small 3};
\node[yshift=6pt]              at (D) {\small 4};
\node[yshift=6pt]              at (E) {\small 5};

\end{tikzpicture}
\end{subfigure}
\hfill
\begin{subfigure}[t]{0.48\textwidth}
\centering
\begin{tikzpicture}[scale=1]

\def\rad{1}

\path[use as bounding box] (-3.15,-2.95) rectangle (3.15,2.95);

\coordinate (O) at (0,0);
\coordinate (A) at (2,0);
\coordinate (B) at ({1},{sqrt(3)});
\coordinate (C) at ({-1},{sqrt(3)});
\coordinate (D) at (-2,0);
\coordinate (E) at ({-1},{-sqrt(3)});

\fill[\cvone]  (O) circle (\rad); 
\fill[\cvtwo]  (A) circle (\rad); 
\fill[\cvthree](B) circle (\rad); 
\fill[\cvfour] (C) circle (\rad); 
\fill[\cvfive] (D) circle (\rad); 
\fill[\cvsix]  (E) circle (\rad); 

\foreach \P in {O,A,B,C,D,E}{
    \draw[circleedge,line width=0.7pt] (\P) circle (\rad);
}

\draw[contactedge,line width=0.8pt] (O)--(A);
\draw[contactedge,line width=0.8pt] (O)--(B);
\draw[contactedge,line width=0.8pt] (O)--(C);
\draw[contactedge,line width=0.8pt] (O)--(D);
\draw[contactedge,line width=0.8pt] (O)--(E);

\draw[contactedge,line width=0.8pt]
    (A)--(B)--(C)--(D)--(E);

\node[xshift=5pt,yshift=-5pt]   at (O) {\small 1};
\node[xshift=6pt,yshift=-3pt]   at (A) {\small 2};
\node[xshift=4pt,yshift=6pt]    at (B) {\small 3};
\node[xshift=-4pt,yshift=6pt]   at (C) {\small 4};
\node[xshift=-6pt]              at (D) {\small 5};
\node[xshift=-4pt,yshift=-6pt]  at (E) {\small 6};

\end{tikzpicture}
\end{subfigure}
\hfill

\caption{The unrestricted global minimisers of the centroidal second moment. Left: the five disk two-row cluster, with \(M=34/5\).
Right: the six disk missing-neighbour flower, with \(M=29/3\).
The segments represent the realised contacts between disk centres.
Both minimisers are unique up to isometry and relabelling.}
\label{fig:all-candidates}
\end{figure}

The global proof reduces the problem to 13 connected outerplanar contact graph classes for five disks and 46 for six disks, counted up to isomorphism. The relevant realisations of these graphs are then analysed while retaining the strict non contact inequalities. For 12 of the five disk graph classes and 44 of the six disk graph classes, validated interval searches combine contact conditions, balance equations, geometric exclusions, and root isolation. The remaining classes are treated by exact algebra and intrinsic second order analysis. Rigorous moment comparisons and a separate argument excluding disconnected equality cases complete the proof, see Section~\ref{sec:verification}.

The relation with rigidity theory is close but distinct. In bar-joint rigidity the graph is prescribed, whereas here it is determined by the active non overlap inequalities of the realised configuration. Classical rigidity and tensegrity ideas therefore enter naturally~\cite{asimowroth,rothwhiteley} together with related ideas from jamming~\cite{donev,connellywhiteley}.

Section~\ref{sec:configuration-space} introduces the hard disk configuration space, realised contact graphs, rolling spaces, and the exact tangent cone. Section~\ref{sec:slice} constructs local Euclidean slices. Section~\ref{sec:first-second} develops centroidal contact balance and the intrinsic Hessian, and Section~\ref{sec:criteria} establishes the leaf obstruction and local minimum criteria. Section~\ref{sec:five} computes the five disk candidate value \(M=34/5\) and proves its strict local minimality. Section~\ref{sec:six} computes the six disk flower value \(M=29/3\) and proves its strict local minimality. Section~\ref{sec:verification} proves the global lower bounds and uniqueness up to isometry and relabelling for both disk numbers, using geometric reductions, exact algebra, and validated interval searches. The Supplementary Notes contain the detailed reductions,
interval coverage, and candidate analysis, while the
repository~\cite{ayalahoffmannvilaxa-dsmp} provides the
computational data, verification code, and reproducibility records.

\section{The hard disk configuration space}
\label{sec:configuration-space}

\subsection{Feasible configurations and compact normalisation}

For \(n\ge1\), define
\begin{equation}
\label{eq:Dn}
 \D_n
 =
 \left\{
 \bc\in(\R^2)^n:
 \|\bc_i-\bc_j\|\ge2
 \text{ for all }i\ne j
 \right\}.
\end{equation}
Thus \(\D_n\) is the full feasible configuration space.

Let
\[
 Z_n
 =
 \left\{
 \bc\in(\R^2)^n:
 \sum_{i=1}^n\bc_i=0
 \right\}.
\]
Every translation orbit in \(\D_n\) meets \(Z_n\) in exactly one point. On
\(Z_n\),
\begin{equation}
\label{eq:M-normalised}
 M(\bc)=\frac12\sum_{i=1}^n\|\bc_i\|^2.
\end{equation}

For \(U>0\), define the normalised compact sublevel
\begin{equation}
\label{eq:compact-sublevel}
 \D_n(U)
 =
 \left\{
 \bc\in\D_n:
 \sum_{i=1}^n\bc_i=0,
 \quad
 M(\bc)\le U
 \right\}.
\end{equation}

\begin{lemma}
\label{lem:existence}
For every \(n\ge1\), the functional \(M\) attains its minimum on \(\D_n\).
Moreover, \(\D_n(U)\) is compact for every \(U>0\).
\end{lemma}

\begin{proof}
The set \(\D_n(U)\) is closed. If \(\bc\in\D_n(U)\), then
\[
 \|\bc_i\|^2
 \le
 \sum_{k=1}^n\|\bc_k\|^2
 =
 2M(\bc)
 \le
 2U
\]
for every \(i\), so the sublevel is bounded. It is therefore compact.

Choose any feasible configuration \(\bc\in\D_n\), and choose
\(U>0\) with \(U\geq M(\bc)\). Every configuration with objective at most
\(U\) is translation equivalent to a point of \(\D_n(U)\). Continuity of
\(M\) and compactness give a minimiser.
\end{proof}

\begin{remark}
The displayed candidates in Figure \ref{fig:all-candidates} give upper bounds \(34/5\) and \(29/3\)
for the five and six disk minima. Hence every global minimiser,
after centroid normalisation, belongs to $\D_5(34/5)$ and $\D_6(29/3)$ respectively. It therefore suffices to exclude competing global
minimisers within these compact sublevel sets.
\end{remark}

\subsection{Exact contacts and realisation sets}

For \(1\le i<j\le n\), define
\begin{equation}
\label{eq:gij}
 g_{ij}(\bc)
 =
 \frac12\|\bc_i-\bc_j\|^2-2.
\end{equation}
For \(i>j\), set \(g_{ij}=g_{ji}\).
Then
\[
 \D_n
 =
 \left\{
 \bc:
 g_{ij}(\bc)\ge0
 \text{ for every }i<j
 \right\},
\]
and \(g_{ij}(\bc)=0\) exactly when disks \(i\) and \(j\) are tangent.

The realised contact graph of \(\bc\in\D_n\) is
\[
 G(\bc)=(V,E(\bc)),
 \qquad
 V=\{1,\ldots,n\},
\]
where
\[
 E(\bc)
 =
 \left\{
 \{i,j\}:
 \|\bc_i-\bc_j\|=2
 \right\}.
\]
When no confusion arises, \(ij\) denotes the unordered edge \(\{i,j\}\). The realised contact graph records the contacting pairs. Its realisation by
the centres \(\bc\) is the contact framework.

For a graph \(G=(V,E)\) on the vertex set
\(V=\{1,\ldots,n\}\), define the exact realisation set
\begin{equation}
\label{eq:exact-realisation}
\mathcal C(G)
 =
 \left\{
 \bc\in\D_n:
 \begin{array}{ll}
 g_{ij}(\bc)=0,
 & \{i,j\}\in E,\\[2mm]
 g_{ij}(\bc)>0,
 & \{i,j\}\notin E
 \end{array}
 \right\}.
\end{equation}

For \(G=(V,E)\), define the active contact map
\begin{equation}
\label{eq:active-contact-map}
 F_G(\bc)
 =
 \bigl(g_{ij}(\bc)\bigr)_{\{i,j\}\in E}.
\end{equation}
When only the edge set is specified, we also write \(F_E\).

The space \(\D_n(U)\) contains configurations of all feasible exact contact
types below the objective level \(U\), whereas \(\mathcal C(G)\) imposes one
fixed exact contact graph. If, for \(\bc\in\mathcal C(G)\), there is a
neighbourhood \(\mathcal U\) of \(\bc\) such that
\(\mathcal C(G)\cap\mathcal U\) is a smooth submanifold of
\((\R^2)^n\), we call this local piece a \emph{realised contact manifold}
through \(\bc\). Note an abstract graph, its exact realisation set, a realised
contact manifold, and a selected representative are distinct objects.

\begin{lemma}
\label{lem:constant-rank-contact-locus}
Let \(G=(V,E)\), and suppose that \(F_G\) has constant rank in a neighbourhood
of \(\bc\in\mathcal C(G)\). Then, after shrinking to a neighbourhood
\(\mathcal U\) of \(\bc\),
\[
 N
 =
 F_G^{-1}(0)\cap\mathcal U
 =
 \mathcal C(G)\cap\mathcal U
\]
is a realised contact manifold through \(\bc\), and
\[
 T_{\bc}N
 =
 \ker DF_G(\bc).
\]
If
\[
 \rank DF_G(\bc)=|E|,
\]
then the active contact constraints satisfy the linear independence constraint
qualification.
\end{lemma}

\begin{proof}
By the constant rank theorem, after restricting to a neighbourhood
\(\mathcal U\) of \(\bc\), the set
\[
 F_G^{-1}(0)\cap\mathcal U
\]
is a smooth submanifold, and its tangent space at \(\bc\) is
\[
 \ker DF_G(\bc).
\]
Since \(\bc\in\mathcal C(G)\), every nonedge \(\{i,j\}\notin E\) has
positive slack at \(\bc\). After shrinking \(\mathcal U\) if necessary,
these strict inequalities remain satisfied throughout \(\mathcal U\).
Therefore
\[
 F_G^{-1}(0)\cap\mathcal U
 =
 \mathcal C(G)\cap\mathcal U
 =
 N.
\]
Thus \(N\) is a realised contact manifold through \(\bc\), with
\[
 T_{\bc}N
 =
 \ker DF_G(\bc).
\]
Full row rank is exactly linear independence of the active contact
differentials.
\end{proof}

\subsection{Rolling spaces and the feasible cone}

For every realised contact \(\{i,j\}\in E(\bc)\), set
\[
 \bu_{ij}
 =
 \frac{\bc_j-\bc_i}{\|\bc_j-\bc_i\|}.
\]

For ambient first order admissibility, contacts may open but may not shorten.
The first order admissible cone is
\begin{equation}
\label{eq:feasible-cone}
 \Adm(\bc)
 =
 \left\{
 \bv\in(\R^2)^n:
 \langle\bu_{ij},\bv_j-\bv_i\rangle\ge0
 \text{ for every }\{i,j\}\in E(\bc)
 \right\}.
\end{equation}
Thus \(\Adm(\bc)\) records the full first order feasibility relation coming
from the non overlap inequalities defining \(\D_n\).

The rolling space is the contact preserving equality part of
\(\Adm(\bc)\):
\begin{equation}
\label{eq:rolling-space}
 \Roll(\bc)
 =
 \left\{
 \bv\in(\R^2)^n:
 \langle\bu_{ij},\bv_j-\bv_i\rangle=0
 \text{ for every }\{i,j\}\in E(\bc)
 \right\}.
\end{equation}
Hence
\[
 \Roll(\bc)\subseteq\Adm(\bc),
\]
and its relation with the admissible cone is
\[
 \Adm(\bc)\cap(-\Adm(\bc))=\Roll(\bc).
\]

Equivalently,
\[
 \Roll(\bc)
 =
 \ker DF_{G(\bc)}(\bc).
\]

By Lemma~\ref{lem:constant-rank-contact-locus}, if \(F_{G(\bc)}\)
has constant rank in a neighbourhood of \(\bc\), then the corresponding
realised contact manifold \(N\) through \(\bc\) satisfies
\[
 T_{\bc}N
 =
 \Roll(\bc).
\]

Let \(J(x,y)=(-y,x)\). The infinitesimal orientation preserving Euclidean
motions are
\begin{equation}
\label{eq:rigid-space}
 R(\bc)
 =
 \left\{
 (a+\omega J\bc_1,\ldots,a+\omega J\bc_n):
 a\in\R^2,\ \omega\in\R
 \right\}.
\end{equation}
Every infinitesimal Euclidean motion preserves all pairwise distances to first
order, and therefore
\[
 R(\bc)\subseteq\Roll(\bc).
\]

Unlike in a bar-joint framework, the contact graph here is realised by \(\bc\) rather than prescribed. We then call infinitesimal rigidity with respect to \(\Roll(\bc)\) rolling rigidity.

\begin{definition}
A configuration \(\bc\in\D_n\) is rolling rigid if
\[
 \Roll(\bc)=R(\bc).
\]
Equivalently, every first order motion preserving all realised contacts is an infinitesimal rigid motion.
\end{definition}

\subsection{The exact tangent cone}

We consider tangent cone
\[
 T_{\D_n}(\bc)
 =
 \left\{
 \bv\in(\R^2)^n:
 \begin{array}{l}
 \text{there exist }t_m\downarrow0
 \text{ and }\bc^{(m)}\in\D_n\\
 \text{such that }
 (\bc^{(m)}-\bc)/t_m\longrightarrow\bv
 \end{array}
 \right\}.
\]

\begin{lemma}
\label{lem:exact-tangent-cone}
For every \(\bc\in\D_n\), the tangent cone is
\[
 T_{\D_n}(\bc)=\Adm(\bc).
\]
\end{lemma}

\begin{proof}
Suppose first that \(\bv\in T_{\D_n}(\bc)\). There are
\(t_m\downarrow0\) and \(\bc^{(m)}\in\D_n\) such that
\[
 \frac{\bc^{(m)}-\bc}{t_m}\longrightarrow\bv.
\]
For every active contact \(\{i,j\}\in E(\bc)\),
\[
 g_{ij}(\bc^{(m)})-g_{ij}(\bc)\ge0.
\]
Divide by \(t_m\) and let \(m\to\infty\). This gives
\[
 Dg_{ij}(\bc)[\bv]\ge0.
\]
At a contact,
\[
 Dg_{ij}(\bc)[\bv]
 =
 \langle\bc_i-\bc_j,\bv_i-\bv_j\rangle
 =
 2\langle\bu_{ij},\bv_j-\bv_i\rangle.
\]
Hence \(\bv\in \Adm(\bc)\).

Conversely, suppose \(\bv\in \Adm(\bc)\). For an active contact, the exact
quadratic expansion is
\begin{equation}
\label{eq:exact-quadratic-contact}
 g_{ij}(\bc+t\bv)
 =
 tDg_{ij}(\bc)[\bv]
 +
 \frac{t^2}{2}\|\bv_i-\bv_j\|^2.
\end{equation}
Both terms on the right are nonnegative for \(t\ge0\). Every inactive pair
has positive slack and remains non overlapping for all sufficiently small
\(t\ge0\). There are only finitely many pairs, so
\[
 \bc+t\bv\in\D_n
\]
for all sufficiently small \(t\ge0\). Thus
\(\bv\in T_{\D_n}(\bc)\).
\end{proof}


\section{Local Euclidean slices}
\label{sec:slice}

A configuration \(\bc\in\D_n\) is an \emph{ambient local minimiser}
of \(M\) if there is a neighbourhood \(\mathcal U\) of \(\bc\)
in \((\R^2)^n\) such that
\(M(\bx)\ge M(\bc)\) for every \(\bx\in\mathcal U\cap\D_n\).

\begin{definition}
\label{def:strict-modulo}
A configuration \(\bc\in\D_n\) is an ambient strict local minimiser of \(M\)
modulo orientation-preserving Euclidean motions if there is a neighbourhood
\(\mathcal U\) of \(\bc\) in \((\R^2)^n\) such that every
\(\bx\in\mathcal U\cap\D_n\) which is not obtained from \(\bc\) by an
orientation-preserving Euclidean motion satisfies
$
M(\bx)>M(\bc).
$
\end{definition}

For \(n\ge2\), the hard disk inequalities imply that the centres are distinct,
and therefore \(R(\bc)\) has dimension three. Let \(Q_\theta\) denote the
counterclockwise rotation through angle \(\theta\).

\begin{lemma}
\label{lem:local-slice}
Let \(n\ge2\) and \(\bc\in\D_n\). With the standard Euclidean inner product
on \((\R^2)^n\), the affine space
\[
 \mathcal S_{\bc}
 =
 \bc+R(\bc)^\perp
\]
is a local slice for the action of orientation preserving Euclidean motions. More precisely, every configuration
\(\bx\) sufficiently close to \(\bc\) has a unique representation
\[
 \bx_i=a+Q_\theta\bs_i,
 \qquad
 \bs\in\mathcal S_{\bc},
\]
with \((\bs,a,\theta)\) sufficiently close to \((\bc,0,0)\).
\end{lemma}

\begin{proof}
Consider the smooth map
\[
 \Phi:
 \mathcal S_{\bc}\times\R^2\times\R
 \longrightarrow
 (\R^2)^n
\]
defined by
\[
 \Phi(\bs,a,\theta)
 =
 (a+Q_\theta\bs_1,\ldots,a+Q_\theta\bs_n).
\]
At \((\bc,0,0)\), its derivative is
\[
 D\Phi(\bc,0,0)[\bz,b,\omega]
 =
 \bz+
 (b+\omega J\bc_1,\ldots,b+\omega J\bc_n).
\]
Since
\[
 (\R^2)^n
 =
 R(\bc)^\perp\oplus R(\bc),
\]
this derivative is an isomorphism. The inverse function theorem
\cite[Chapter~4]{lee} gives the asserted local diffeomorphism and uniqueness.
\end{proof}

\begin{lemma}
\label{lem:slice-reduction}
A configuration \(\bc\in\D_n\), \(n\ge2\), is an ambient strict local
minimiser of \(M\) modulo orientation preserving Euclidean motions if and only if it is a strict local minimiser
of \(M\) on
\[
 \D_n\cap\mathcal S_{\bc}.
\]
\end{lemma}

\begin{proof}
Every nearby configuration has a unique representation
\[
 \bx=\Phi(\bs,a,\theta)
\]
from Lemma~\ref{lem:local-slice}. The feasible set and the objective are
invariant under orientation preserving Euclidean motions, so
\[
 \bx\in\D_n
 \quad\Longleftrightarrow\quad
 \bs\in\D_n
\]
and
\[
 M(\bx)=M(\bs).
\]
By local uniqueness, \(\bx\) lies in the orbit of \(\bc\) if and only if
\(\bs=\bc\). The two strict local minimum statements are therefore equivalent.
\end{proof}

\section{First variation, contact balance, and intrinsic Hessian}
\label{sec:first-second}

\subsection{First variation of the centroidal second moment}

For \(\bv=(\bv_1,\ldots,\bv_n)\), its centroid is
$$\vbar=\frac1n\sum_{i=1}^n\bv_i.$$ Along \(\bc+t\bv\), the centroid is
$ \cbar(t)=\cbar+t\vbar$. Differentiating \eqref{eq:objective} gives
\begin{equation}
\label{eq:first-variation}
 DM(\bc)[\bv]
 =
 \sum_{i=1}^n
 \langle\bc_i-\cbar,\bv_i-\vbar\rangle.
\end{equation}
Since
\[
 \sum_{i=1}^n(\bc_i-\cbar)=0,
\]
this simplifies to
\begin{equation}
\label{eq:first-variation-simple}
 DM(\bc)[\bv]
 =
 \sum_{i=1}^n
 \langle\bc_i-\cbar,\bv_i\rangle.
\end{equation}
Thus
\begin{equation}
\label{eq:gradient}
 \nabla_{\bc_i}M(\bc)=\bc_i-\cbar.
\end{equation}

For \(\bw=(\bw_1,\ldots,\bw_n)\), set
\(\wbar=\frac1n\sum_{i=1}^n\bw_i\).
The ambient Hessian is constant:
\begin{equation}
\label{eq:ambient-hessian}
 D^2M(\bc)[\bv,\bw]
 =
 \sum_{i=1}^n
 \langle\bv_i-\vbar,\bw_i-\wbar\rangle.
\end{equation}
In particular,
\begin{equation}
\label{eq:ambient-hessian-diagonal}
 D^2M(\bc)[\bv,\bv]
 =
 \sum_{i=1}^n
 \|\bv_i-\vbar\|^2
 =
 \frac1n\sum_{i<j}\|\bv_i-\bv_j\|^2.
\end{equation}

\subsection{Centroidal contact balance}

Suppose \(F_E\) has constant rank in an ambient neighbourhood of
\(\bc\), where \(E\subseteq E(\bc)\), and suppose that \(\bc\) is
critical for the restriction of \(M\) to the resulting local equality
manifold \(F_E^{-1}(0)\). Then the differential of \(M\) belongs to the span of
the active contact differentials. Hence there are real multipliers
\(\lambda_{ij}\), not necessarily nonnegative, such that
\begin{equation}
\label{eq:balance-vector}
 \nabla M(\bc)
 =
 \sum_{\{i,j\}\in E}
 \lambda_{ij}\nabla g_{ij}(\bc).
\end{equation}
At each vertex this is
\begin{equation}
\label{eq:balance-vertex}
 \bc_i-\cbar
 =
 \sum_{j:\{i,j\}\in E}
 \lambda_{ij}(\bc_i-\bc_j),
 \qquad
 \lambda_{ij}=\lambda_{ji}.
\end{equation}
We call \eqref{eq:balance-vertex} the \emph{centroidal contact balance}.

For ambient minimality the coefficients can be chosen nonnegative, and this
conclusion does not require a nonlinear constraint qualification.

\begin{proposition}
\label{prop:ambient-first-order}
If \(\bc\) is an ambient local minimiser of \(M\) on \(\D_n\), then
\begin{equation}
\label{eq:ambient-first-order}
 DM(\bc)[\bv]\ge0
 \qquad
 \text{for every }\bv\in \Adm(\bc).
\end{equation}
Equivalently, there are multipliers
\[
 \lambda_{ij}\ge0,
 \qquad
 \{i,j\}\in E(\bc),
\]
satisfying the centroidal contact balance
\eqref{eq:balance-vertex}.
\end{proposition}

\begin{proof}
By Lemma~\ref{lem:exact-tangent-cone}, every \(\bv\in \Adm(\bc)\) generates a
feasible ray \(\bc+t\bv\) for sufficiently small \(t\ge0\). Local minimality
therefore gives \eqref{eq:ambient-first-order}.

The cone \(\Adm(\bc)\) is polyhedral. The finite dimensional polar cone identity,
equivalently Farkas' lemma, gives
\[
 \left\{
 \ell:
 \ell(\bv)\ge0
 \text{ for every }\bv\in \Adm(\bc)
 \right\}
 =
 \operatorname{cone}
 \left\{
 Dg_{ij}(\bc):
 \{i,j\}\in E(\bc)
 \right\};
\]
see, for example, \cite[Chapter~6]{rockafellarwets}. Applying this identity to
\(DM(\bc)\) gives nonnegative multipliers and
\eqref{eq:balance-vertex}.
\end{proof}

\subsection{The intrinsic Hessian on a realised contact manifold}

The full active contact system may be redundant even when its equality set is
smooth. Therefore the intrinsic second derivative must be defined from the
restricted objective itself, not from an arbitrary multiplier representation.

\begin{definition}
\label{def:intrinsic-hessian}
Let \(N\subset(\R^2)^n\) be a smooth submanifold, let \(\bc\in N\), and suppose
\[
 d(M|_N)(\bc)=0.
\]
Choose a smooth local parametrisation
\[
 \phi:U\subset\R^d\longrightarrow N,
 \qquad
 \phi(0)=\bc,
\]
whose derivative at \(0\) is an isomorphism onto \(T_{\bc}N\). For
\[
 \bv=D\phi(0)\xi,
 \qquad
 \bw=D\phi(0)\eta,
\]
define
\begin{equation}
\label{eq:intrinsic-hessian-definition}
 \Hess_N M(\bc)[\bv,\bw]
 =
 D^2(M\circ\phi)(0)[\xi,\eta].
\end{equation}
\end{definition}

\begin{lemma}
\label{lem:hessian-coordinate-independent}
The bilinear form in Definition~\ref{def:intrinsic-hessian} is independent of
the chosen local parametrisation.
\end{lemma}

\begin{proof}
Let \(\psi\) be another local parametrisation. The transition map between the
two coordinate systems is a local diffeomorphism. Under a change of variables,
the second derivative acquires an additional term involving the first
derivative of \(M|_N\). That first derivative vanishes at \(\bc\), so the
additional term is zero. The two coordinate expressions therefore define the
same bilinear form on \(T_{\bc}N\).
\end{proof}

\begin{proposition}
\label{prop:independent-defining-equations}
Suppose \(N\subset(\R^2)^n\) is locally defined near \(\bc\) by independent smooth
equations
\[
 h_1=\cdots=h_r=0
\]
and
\[
 dM(\bc)=\sum_{a=1}^r\mu_a\,dh_a(\bc)
\]
as covectors on \((\R^2)^n\). Then, for \(\bv\in T_{\bc}N\),
\begin{equation}
\label{eq:intrinsic-hessian-lagrangian}
 \Hess_N M(\bc)[\bv,\bv]
 =
 D^2M(\bc)[\bv,\bv]
 -
 \sum_{a=1}^r
 \mu_aD^2h_a(\bc)[\bv,\bv].
\end{equation}
\end{proposition}

\begin{proof}
Let \(\gamma:(-\varepsilon,\varepsilon)\to N\) be a smooth curve with
\[
 \gamma(0)=\bc,
 \qquad
 \gamma'(0)=\bv.
\]
Since \(h_a(\gamma(t))=0\),
\[
 dh_a(\bc)[\gamma''(0)]
 =
 -
 D^2h_a(\bc)[\bv,\bv].
\]
Differentiating \(M(\gamma(t))\) twice gives
\[
 \frac{d^2}{dt^2}M(\gamma(t))\bigg|_{t=0}
 =
 D^2M(\bc)[\bv,\bv]
 +
 dM(\bc)[\gamma''(0)].
\]
Substitute
\[
 dM(\bc)=\sum_{a=1}^r\mu_a\,dh_a(\bc)
\]
and use the preceding identity. This gives
\eqref{eq:intrinsic-hessian-lagrangian}.
\end{proof}

\begin{remark}[Redundant active systems]
Where the active contact map has constant rank in an ambient neighbourhood,
one may choose an independent local subfamily of the active contact equations and use
Proposition~\ref{prop:independent-defining-equations}. If the complete active
contact system is redundant, an arbitrary multiplier vector for the full
system must not be inserted into a Lagrangian Hessian without proving
multiplier independence. The validated computation therefore uses either:

\begin{enumerate}[label=\textup{(\roman*)}]
\item a validated local parametrisation of the actual realised contact manifold; or
\item an independent local defining subsystem.
\end{enumerate}

This convention applies to regular redundant systems as well as to singular
rank pieces. On a singular rank piece, option~(ii) requires verification that
the chosen equations locally define the actual realised contact manifold.
\end{remark}

Suppose that \(N\) is locally defined near a critical point \(\bc\) of
\(M|_N\) by an independent contact subsystem \(E_0\). Choose the multipliers
\(\lambda_{ij}\) so that \eqref{eq:balance-vector} holds with \(E=E_0\).
For \(\bv\in T_{\bc}N\), the intrinsic Hessian satisfies
\begin{equation}
\label{eq:contact-intrinsic-hessian}
 \Hess_N M(\bc)[\bv,\bv]
 =
 \sum_{i=1}^n\|\bv_i-\vbar\|^2
 -
 \sum_{\{i,j\}\in E_0}
 \lambda_{ij}\|\bv_i-\bv_j\|^2.
\end{equation}
where \(E_0\) is the chosen independent local defining subsystem.

\begin{lemma}
\label{lem:second-order-exclusion}
Let \(N\) be a realised contact manifold through \(\bc\), and suppose that
\[
 d(M|_N)(\bc)=0.
\]
If there is a tangent vector \(\bv\in T_{\bc}N\) such that
\[
 \Hess_N M(\bc)[\bv,\bv]<0,
\]
then \(\bc\) is not a local minimiser of \(M\) on \(N\), and hence is not a
global minimiser on \(\D_n\).
\end{lemma}

\begin{proof}
Choose a local chart for \(N\) and represent \(\bv\) by a coordinate vector.
The ordinary second derivative test gives a smooth curve in \(N\) through
\(\bc\) along which \(M\) decreases.
\end{proof}

\section{Contact leaf obstruction and local minimum criteria}
\label{sec:criteria}

\subsection{The ambient oriented leaf obstruction}

A vertex \(i\) is a \emph{contact leaf} if it has exactly one contact
neighbour.

\begin{lemma}
\label{lem:ambient-leaf}
Let \(\bc\in\D_n\), and suppose that \(i\) is a contact leaf with unique
neighbour \(j\). If \(\bc\) is an ambient local minimiser of \(M\), then there
is \(\alpha\ge0\) such that
\begin{equation}
\label{eq:leaf-ray}
 \bc_i-\cbar
 =
 \alpha(\bc_i-\bc_j).
\end{equation}
Equivalently, the centroid lies on the ray emanating from the leaf centre in
the direction of its contact neighbour.

If \eqref{eq:leaf-ray} fails, there is a feasible first order direction
\(\bv\in \Adm(\bc)\) with
\[
 DM(\bc)[\bv]<0.
\]
\end{lemma}

\begin{proof}
Set
\[
 p=\bc_i-\cbar,
 \qquad
 a=\bc_i-\bc_j.
\]
For \(w\in\R^2\), define a velocity supported at the leaf by
\[
 \bv_i=w,
 \qquad
 \bv_k=0
 \quad(k\ne i).
\]
The only active inequality involving \(i\) gives
\[
 Dg_{ij}(\bc)[\bv]
 =
 \langle a,w\rangle.
\]
Thus \(\bv\in \Adm(\bc)\) exactly when
\[
 \langle a,w\rangle\ge0.
\]
The first variation is
\[
 DM(\bc)[\bv]=\langle p,w\rangle.
\]

The condition
\[
 \langle p,w\rangle\ge0
 \quad
 \text{whenever}
 \quad
 \langle a,w\rangle\ge0
\]
holds if and only if
\[
 p\in\R_{\ge0}a.
\]
This is precisely \eqref{eq:leaf-ray}. If it fails, separation of \(p\) from
the ray \(\R_{\ge0}a\) gives an admissible \(w\) with
\[
 \langle a,w\rangle\ge0,
 \qquad
 \langle p,w\rangle<0.
\]
\end{proof}

\begin{lemma}
\label{lem:leaf-second-order-bound}
Let \(\bc\in\D_n\) be an ambient local minimiser of \(M\), and suppose that
\(i\) is a contact leaf with unique neighbour \(j\). Write
\[
 \bc_i-\cbar
 =
 \alpha(\bc_i-\bc_j),
 \qquad
 \alpha\ge0,
\]
as in Lemma~\ref{lem:ambient-leaf}. Then
\[
 0\le\alpha\le\frac{n-1}{n}.
\]
\end{lemma}

\begin{proof}
Set
\[
 a=\bc_i-\bc_j,
\]
so that \(\|a\|=2\), and rotate only the leaf centre about its contact
neighbour:
\[
 \bc_i(t)=\bc_j+Q_ta,
 \qquad
 \bc_k(t)=\bc_k
 \quad(k\ne i).
\]
For sufficiently small \(|t|\), this is a feasible curve. The unique active
contact incident to \(i\) is preserved exactly, every other active contact is
unchanged, and all inactive pairs retain positive slack.

At \(t=0\),
\[
 \bc_i'(0)=Ja,
 \qquad
 \bc_i''(0)=-a.
\]
The first derivative vanishes:
\[
 \frac{d}{dt}M(\bc(t))\bigg|_{t=0}
 =
 \langle\bc_i-\cbar,Ja\rangle
 =
 \alpha\langle a,Ja\rangle
 =
 0.
\]
Using \eqref{eq:ambient-hessian-diagonal},
\[
 D^2M(\bc)[\bc'(0),\bc'(0)]
 =
 \frac{n-1}{n}\|a\|^2,
\]
while
\[
 DM(\bc)[\bc''(0)]
 =
 -\alpha\|a\|^2.
\]
Hence
\[
 \frac{d^2}{dt^2}M(\bc(t))\bigg|_{t=0}
 =
 \|a\|^2
 \left(
 \frac{n-1}{n}-\alpha
 \right).
\]
Ambient local minimality forces this quantity to be nonnegative, and
therefore
\[
 \alpha\le\frac{n-1}{n}.
\]
\end{proof}

\begin{remark}
Comparing \eqref{eq:leaf-ray} with the centroidal contact balance
\eqref{eq:balance-vertex}, because vertex \(i\) has degree one, the sum
collapses to a single term and the scalar \(\alpha\) is identically the
equilibrium contact multiplier \(\lambda_{ij}\) on the unique leaf edge.
Lemma~\ref{lem:leaf-second-order-bound} strengthens the first-order ray
condition at an ambient local minimiser: since
\[
 0\le\alpha\le\frac{n-1}{n}<1,
\]
the centroid lies on the closed segment from the leaf centre toward its
contact neighbour, strictly before the neighbour.
\end{remark}

\subsection{The trivial critical cone criterion}

Define the critical cone
\begin{equation}
\label{eq:critical-cone}
\Adm_M(\bc)
=
\left\{
\bv\in\Adm(\bc):
DM(\bc)[\bv]=0
\right\}.
\end{equation}

\begin{theorem}
\label{thm:critical-cone}
Let \(\bc\in\D_n\), \(n\ge2\). Suppose
\begin{equation}
\label{eq:first-order-nonnegative}
 DM(\bc)[\bv]\ge0
 \qquad
 \text{for every }\bv\in \Adm(\bc),
\end{equation}
and
\begin{equation}
\label{eq:trivial-reduced-cone}
 \Adm_M(\bc)\cap R(\bc)^\perp=\{0\}.
\end{equation}
Then \(\bc\) is an ambient strict local minimiser of \(M\)
modulo orientation-preserving Euclidean motions.
No constraint qualification is required.
\end{theorem}

\begin{proof}
Pass to the slice
\[
 \mathcal S_{\bc}=\bc+R(\bc)^\perp.
\]
Suppose that strict local minimality fails on the slice. Then there are
points \(\bc^{(m)}\ne\bc\)
\[
 \bc^{(m)}\in\D_n\cap\mathcal S_{\bc}
\]
such that
\[
 \bc^{(m)}\longrightarrow\bc,
 \qquad
 M(\bc^{(m)})\le M(\bc).
\]
Set
\[
 t_m=\|\bc^{(m)}-\bc\|,
 \qquad
 \bw^{(m)}
 =
 \frac{\bc^{(m)}-\bc}{t_m}.
\]
After taking a subsequence,
\[
 \bw^{(m)}\longrightarrow\bw,
 \qquad
 \|\bw\|=1.
\]
Because the sequence lies in the affine slice,
\[
 \bw\in R(\bc)^\perp.
\]
Feasibility gives
\[
 \bw\in T_{\D_n}(\bc)=\Adm(\bc).
\]

Taylor expansion gives
\[
 M(\bc^{(m)})
 =
 M(\bc)
 +
 t_mDM(\bc)[\bw^{(m)}]
 +
 O(t_m^2).
\]
Since \(M(\bc^{(m)})\le M(\bc)\), divide by \(t_m\) and pass to the limit:
\[
 DM(\bc)[\bw]\le0.
\]
The first order condition \eqref{eq:first-order-nonnegative} gives equality.
Therefore
\[
 \bw\in \Adm_M(\bc)\cap R(\bc)^\perp,
\]
contradicting \eqref{eq:trivial-reduced-cone} and \(\|\bw\|=1\).

Thus \(\bc\) is a strict local minimiser on the slice. The conclusion modulo orientation preserving Euclidean motions follows from Lemma~\ref{lem:slice-reduction}.
\end{proof}

\begin{corollary}
\label{cor:rigid-positive-balance}
Let \(\bc\in\D_n\), \(n\ge2\), be rolling rigid. Suppose that there are
strictly positive multipliers
\[
 \lambda_{ij}>0,
 \qquad
 \{i,j\}\in E(\bc),
\]
satisfying the centroidal contact balance
\[
 \bc_i-\cbar
 =
 \sum_{j:\{i,j\}\in E(\bc)}
 \lambda_{ij}(\bc_i-\bc_j)
\]
at every vertex. Then \(\bc\) is an ambient strict local minimiser of \(M\)
modulo orientation preserving Euclidean motions.
\end{corollary}

\begin{proof}
For every \(\bv\in \Adm(\bc)\), the balance equation gives
\[
 DM(\bc)[\bv]
 =
 \sum_{\{i,j\}\in E(\bc)}
 \lambda_{ij}Dg_{ij}(\bc)[\bv]
 \ge0.
\]
If equality holds, strict positivity of every multiplier forces
\[
 Dg_{ij}(\bc)[\bv]=0
 \qquad
 \text{for every }\{i,j\}\in E(\bc).
\]
Hence
\[
 \bv\in\Roll(\bc)=R(\bc).
\]
Therefore
\[
 \Adm_M(\bc)\cap R(\bc)^\perp=\{0\}.
\]
Apply Theorem~\ref{thm:critical-cone}.
\end{proof}

\subsection{Rigidity under two contact attachment}

\begin{lemma}
\label{lem:triangle-rigid}
Let three centres form a noncollinear triangle whose three edges are realised
contacts. Then the contact framework is rolling rigid.
\end{lemma}

\begin{proof}
There are three contact equations in a six dimensional velocity space. Their
gradients are independent. Indeed, if a linear combination of the three edge
gradients vanishes, then at one vertex the two incident edge vectors occur in
a linear combination equal to zero. Since those vectors are independent, the
two corresponding coefficients vanish; the third then vanishes as well.

The contact Jacobian therefore has rank three, and its kernel has dimension
three. Since the three dimensional rigid motion space is contained in the
kernel, the two spaces are equal.
\end{proof}

\begin{lemma}
\label{lem:two-contact-attachment}
Let a realised contact framework on centres
\[
 \bc_1,\ldots,\bc_m
\]
be rolling rigid. Add a new centre \(\bc_q\) with contacts to two old centres
\(\bc_i,\bc_j\). If
\[
 \bc_q-\bc_i
 \qquad\text{and}\qquad
 \bc_q-\bc_j
\]
are linearly independent, then the enlarged framework is rolling rigid.
Additional realised contacts preserve this conclusion.
\end{lemma}

\begin{proof}
Let \(\bv\) be a rolling velocity of the enlarged framework. Its restriction
to the old vertices is an infinitesimal Euclidean motion, so there are
\(a\in\R^2\) and \(\omega\in\R\) such that
\[
 \bv_k=a+\omega J\bc_k
 \qquad
 (1\le k\le m).
\]
Set
\[
 w=\bv_q-(a+\omega J\bc_q).
\]
Subtracting the same infinitesimal Euclidean motion from all velocities, the
new rolling equations become
\[
 \langle\bc_q-\bc_i,w\rangle=0,
 \qquad
 \langle\bc_q-\bc_j,w\rangle=0.
\]
The two contact directions are independent, so \(w=0\). Thus the new velocity
also agrees with the same infinitesimal Euclidean motion. Additional contacts
only add rolling equations and cannot enlarge the rolling space.
\end{proof}

\begin{remark}
Two contact attachment propagates rigidity only. Adding a centre changes the
centroidal load, so the contact balance equations must be solved again for the
enlarged configuration.
\end{remark}

\section{The five disk solution of Graham and Sloane and its local minimality}
\label{sec:five}
We begin with the five disk configuration of Graham and Sloane as an
application of the local theory. This section computes its objective value
and proves strict local minimality. The unrestricted global lower bound and
uniqueness are established in Section~\ref{sec:verification} through the
census of realised contact graphs and the subsequent exact analysis.

Consider the centroid normalised configuration
\begin{equation}
\label{eq:five-coordinates}
\begin{aligned}
 \bc_1&=(-2,-2\sqrt3/5),
&
 \bc_2&=(0,-2\sqrt3/5),
&
 \bc_3&=(2,-2\sqrt3/5),
\\
 \bc_4&=(-1,3\sqrt3/5),
&
 \bc_5&=(1,3\sqrt3/5).
\end{aligned}
\end{equation}
Its contact edge set is
\begin{equation}
\label{eq:five-edges}
 E_5=\{12,23,45,14,24,25,35\}.
\end{equation}
The three non edges satisfy
\[
 \|\bc_1-\bc_3\|=4,
 \qquad
 \|\bc_1-\bc_5\|
 =
 \|\bc_3-\bc_4\|
 =
 2\sqrt3.
\]
Thus \eqref{eq:five-edges} is the exact realised contact edge set.

Since the centroid is zero,
\[
 \sum_{i=1}^5\|\bc_i\|^2=\frac{68}{5},
\]
and therefore
\begin{equation}
\label{eq:five-value}
 M(\bc)=\frac{34}{5}.
\end{equation}

\begin{proposition}
\label{prop:five-local-data}
The five disk contact framework is rolling rigid. Its seven active contact
gradients are linearly independent. The centroidal contact balance is solved by
\[
 \lambda_{12}=\lambda_{23}=\frac45,
 \qquad
 \lambda_{45}=\frac35,
\]
and
\[
 \lambda_{14}=\lambda_{35}=\frac25,
 \qquad
 \lambda_{24}=\lambda_{25}=\frac15.
\]
All multipliers are strictly positive.
\end{proposition}

\begin{proof}
The vertices \(1,2,4\) form a noncollinear contact triangle and are rolling
rigid by Lemma~\ref{lem:triangle-rigid}. Vertex \(5\) is attached to vertices
\(2,4\) in two independent contact directions. Vertex \(3\) is then attached
to vertices \(2,5\) in two independent directions. Two applications of
Lemma~\ref{lem:two-contact-attachment} prove
\[
 \Roll(\bc)=R(\bc).
\]

The ambient velocity space has dimension ten, while \(R(\bc)\) has dimension
three. Hence the contact Jacobian has rank seven, equal to the number of active
contacts, so the active gradients are linearly independent.

Substitution of the displayed multipliers and the coordinates
\eqref{eq:five-coordinates} into
\[
 \bc_i
 =
 \sum_{j:\{i,j\}\in E_5}
 \lambda_{ij}(\bc_i-\bc_j)
\]
verifies the five vertex equations exactly.
\end{proof}

\begin{corollary}
\label{cor:five-local}
The configuration \eqref{eq:five-coordinates} is an ambient strict local
minimiser of \(M\) modulo orientation preserving Euclidean motions.
\end{corollary}

\begin{proof}
Apply Corollary~\ref{cor:rigid-positive-balance}.
\end{proof}

\section{The six disk flower and its local minimality}
\label{sec:six}
We now compute \(M=29/3\) for the six disk flower and prove its strict local
minimality. Section~\ref{sec:verification} establishes the global lower bound
and uniqueness..

Consider
\begin{equation}
\label{eq:six-coordinates}
\begin{aligned}
 \bc_1&=(0,0),
&
 \bc_2&=(2,0),
&
 \bc_3&=(1,\sqrt3),
\\
 \bc_4&=(-1,\sqrt3),
&
 \bc_5&=(-2,0),
&
 \bc_6&=(-1,-\sqrt3).
\end{aligned}
\end{equation}
Its centroid is
\[
 \cbar
 =
 \left(-\frac16,\frac{\sqrt3}{6}\right).
\]
Its exact contact edge set is
\begin{equation}
\label{eq:six-edges}
 E_6
 =
 \{12,13,14,15,16,23,34,45,56\}.
\end{equation}
Every non edge distance is \(2\sqrt3\) or \(4\), and hence is strictly larger
than \(2\).

Moreover,
\[
 \sum_{i=1}^6\|\bc_i\|^2=20,
 \qquad
 \|\cbar\|^2=\frac19.
\]
Therefore
\begin{equation}
\label{eq:six-value}
 M(\bc)
 =
 \frac12
 \left(
 20-6\cdot\frac19
 \right)
 =
 \frac{29}{3}.
\end{equation}

\begin{proposition}
\label{prop:six-local-data}
The six disk contact framework is rolling rigid. Its nine active contact
gradients are linearly independent. The centroidal contact balance is solved by
\[
 \lambda_{12}=1,
 \quad
 \lambda_{13}=\frac23,
 \quad
 \lambda_{14}=\frac12,
 \quad
 \lambda_{15}=\frac23,
 \quad
 \lambda_{16}=1,
\]
and
\[
 \lambda_{23}=\frac16,
 \quad
 \lambda_{34}=\frac13,
 \quad
 \lambda_{45}=\frac13,
 \quad
 \lambda_{56}=\frac16.
\]
All multipliers are strictly positive.
\end{proposition}

\begin{proof}
The vertices \(1,2,3\) form a noncollinear contact triangle. Attach vertex
\(4\) to vertices \(1,3\), then vertex \(5\) to vertices \(1,4\), and finally
vertex \(6\) to vertices \(1,5\). At every step the two new contact directions
are independent. Lemmas~\ref{lem:triangle-rigid} and
\ref{lem:two-contact-attachment} give
\[
 \Roll(\bc)=R(\bc).
\]

The ambient velocity space has dimension twelve and the rolling space has
dimension three. Thus
\[
 \rank DF_{E_6}(\bc)=9,
\]
equal to the number of active contacts, so the active gradients are linearly
independent.

Substitution of the displayed multipliers into
\[
 \bc_i-\cbar
 =
 \sum_{j:\{i,j\}\in E_6}
 \lambda_{ij}(\bc_i-\bc_j)
\]
verifies every vertex equation exactly.
\end{proof}

\begin{corollary}
\label{cor:six-local}
The configuration \eqref{eq:six-coordinates} is an ambient strict local
minimiser of \(M\) modulo orientation preserving Euclidean motions.
\end{corollary}

\begin{proof}
Apply Corollary~\ref{cor:rigid-positive-balance}.
\end{proof}

\section{Global classification through contact graphs}
\label{sec:verification}
The five and six disk candidate values of \(M\) are \(34/5\) and
\(29/3\), respectively. We denote these comparison values by \(U_5\)
and \(U_6\). We now prove the global lower bounds and classify the equality cases. The geometric reductions below give a
finite list of contact graphs. Algebraic arguments and interval searches
then account for every possible ambient local minimum on that list.
The Supplementary Notes contain the detailed search arguments,
candidate analyses, and computational record.

\subsection{Structural properties of minimisers}
\label{subsec:structura-properties-minimizer}

A global minimiser exists by Lemma~\ref{lem:existence}, and its second moment
is at most \(U_n\). We now show that it suffices to consider packings
with connected, outerplanar contact graphs.
  
\begin{lemma}[Reduction to connected contact graphs]
\label{lem:connected-census-reduction}
For \(n\in\{5,6\}\), a disconnected global minimiser can be moved,
without changing its second moment, to a packing with fewer contact
components. Consequently, a connected global minimiser exists.
\end{lemma}

\begin{proof}
Let \(\bc\) be a disconnected global minimiser, so \(M(\bc)\le U_n\).
Small translations of each contact component are feasible in every
direction, since different components have no contacts. By the gradient
formula \eqref{eq:gradient}, stationarity under these translations
forces every component centroid to equal \(\cbar\). Put their common
centroid at the origin. Rotating one component about the origin
preserves its internal contacts, the total centroid, and the second moment.

Choose a centre \(p\) of smallest radius \(r\), and let \(C\) be its
component. Suppose that rotating \(C\) through a full turn creates no
new contact. By continuity, all distances between \(C\) and the other
components remain greater than two. For an outside centre of radius
\(s\ge r\), its minimum distance from the rotating centre \(p\) is
\(s-r\). Thus every outside centre has radius greater than \(r+2\).
There are at least two outside centres: otherwise the sole outside
centre would form a singleton component at the origin, contradicting
this inequality.

Every centre other than \(p\), of radius \(s\), satisfies
\(s\ge\max\{r,2-r\}\ge|2-r|\), by minimality of \(r\) and separation
from \(p\). At least two of these centres lie outside \(C\), so their radii exceed \(r+2\). Hence,
\[
\begin{aligned}
M
&>\frac12\left[r^2+(n-3)(2-r)^2+2(r+2)^2\right]\\
&=
\begin{cases}
8+\dfrac52r^2, & n=5,\\[2mm]
\dfrac{29}{3}+3\left(r-\dfrac13\right)^2, & n=6.
\end{cases}
\end{aligned}
\]
In either case this contradicts \(M\le U_n\).

A new contact must therefore occur. Stop at the first one. All
previous contacts are preserved, so the resulting packing has fewer
contact components and is again a global minimiser. Repetition gives
a connected global minimiser.
\end{proof}

\begin{lemma}[Reduction to outerplanar graphs]
\label{lem:outerplanar-census-reduction}
Every contact graph of at most six disks is outerplanar.
\end{lemma}

\begin{proof}
Draw contact edges as straight segments. No centre lies in the
interior of an edge, since its distance to either endpoint would
be less than two. This also excludes collinear overlap. If two
independent edges \(ab\) and \(cd\) cross, the strict triangle
inequalities through their intersection give
\[
 \|a-c\|+\|b-d\|<\|a-b\|+\|c-d\|=4,
\]
contradicting disk separation. Thus the drawing is a plane
embedding.

Suppose that a centre \(x\) is not on the boundary of the
unbounded face. Then a simple contact cycle
\(v_1,\ldots,v_m\) encloses \(x\), with indices read cyclically.
In each triangle \(xv_iv_{i+1}\), the side \(v_iv_{i+1}\) has
length two and the other sides have length at least two.
The angle at \(x\) is therefore at most \(\pi/3\); in the
collinear case it is zero.

Let \(\Delta_i\) be the change of argument about \(x\) along
the directed edge \(v_iv_{i+1}\). Then
\(|\Delta_i|\le\pi/3\). Since the cycle winds once around \(x\),
\[
 2\pi
 =\left|\sum_{i=1}^m\Delta_i\right|
 \le\sum_{i=1}^m|\Delta_i|
 \le\frac{m\pi}{3}.
\]
Hence \(m\ge6\), requiring at least seven centres including
\(x\), a contradiction. Therefore every vertex is on the
boundary of the unbounded face.
\end{proof}

\subsection{Finite graph census}
\label{subsec:census}

The preceding lemmas reduce the search to connected outerplanar contact
graphs. We used \texttt{plantri} to generate connected planar graphs
on five and six vertices and \texttt{shortg} to remove isomorphic
duplicates \cite{brinkmannmckay,mckaypiperno}. The resulting lists
contain 20 and 99 graphs, respectively. The outerplanarity filter
leaves 13 and 46 (Supplementary Notes, Section~2). Contact equations
and strict noncontact inequalities are imposed in the subsequent
geometric analysis.

Table~\ref{tab:proof-pipeline-summary} divides the graphs by their
verification route. The target graphs are verified algebraically, and
one six-disk graph is analysed through its critical family. All remaining
graphs enter the same interval search. It combines geometric and
contact-balance tests with additional conditions for leaves and bridges
when applicable, as described in Section~\ref{subsec:phase1-search}.

\begin{table}[htbp]
\centering\small
\caption{Partition of connected outerplanar contact-graph classes by
verification route. Every entry is a number of graph classes.}
\label{tab:proof-pipeline-summary}
\begin{tabular}{@{}lrr@{}}
\toprule
Route & \(n=5\) & \(n=6\)\\
\midrule
Interval search & 12 & 44\\
Target graph verified algebraically & 1 & 1\\
Critical family & 0 & 1\\
\midrule
Total connected classes & 13 & 46\\
\bottomrule
\end{tabular}
\end{table}

\subsection{Interval search and coverage}
\label{subsec:phase1-search}

The interval search tests necessary conditions for local minimality
of the second moment. It excludes angular boxes or retains isolated
roots for candidate analysis,
with and without leaves.

For a contact graph \(G=(V,E)\), choose a rooted spanning tree and
fix the centres of its first contact at \((0,0)\) and \((2,0)\).
For a graph with leaves, the first edge is incident to the lowest-labelled
leaf. Otherwise it is the first stored contact edge. Complete the tree
breadth first in label order. For each remaining parent--child edge, its direction \(\theta\) places the next centre:
\[
 \bc_j=\bc_i+2(\cos\theta,\sin\theta).
\]
For \(n\) disks, this leaves \(n-2\) free angles and keeps every pair joined by a tree edge at distance \(2\). Choosing a value for every free angle determines all
centre positions.

We group these choices by assigning an interval to each free angle.
An \emph{angular box} contains all combinations of values from these
intervals. The initial intervals cover a full turn for every free
angle except one, restricted to a half-turn by reflection in the fixed
edge. The resulting box covers every realisation up to isometry.
For graphs with leaves, further graph symmetries reduce the domain
while preserving coverage (Supplementary Notes, Section~4.4).

Interval arithmetic bounds expressions for all combinations in a box
at once, using 512-bit \texttt{Arb} arithmetic with outward rounding
\cite{johansson}. These bounds can exclude whole boxes. Subdivision
can sharpen inconclusive bounds while preserving coverage.

For each box, evaluate \(g_{ij}\) by interval arithmetic for every contact outside the tree and every nonedge. Contacts require \(g_{ij}=0\), whereas nonedges require \(g_{ij} > 0\). Discard the box if a contact interval excludes
zero or a separation interval has no positive values. Each discard
ends the examination of that box.

When the graph has leaves, test their balance and the necessary bounds
\(0\le\lambda_{ij}\le(n-1)/n\)
(Lemmas~\ref{lem:ambient-leaf} and~\ref{lem:leaf-second-order-bound}).
For each bridge, test equilibrium across it and second variation under
rotation of one component. These tests use the centre coordinates directly.
When applicable, also test the bridge equations using the Krawczyk
exclusion criterion described below
(Supplementary Notes, Sections~4 and~6.2).

Surviving boxes are then tested against the contact-balance equations
from Proposition~\ref{prop:ambient-first-order}:
\[
 DF_G(\bc)^T\blambda=\nabla M(\bc),
 \qquad \blambda\geq0.
\]
Apply interval Gaussian elimination to all coordinate equations,
including every contact multiplier. Each pivot must exclude zero.
Failure to find one causes box subdivision. A remaining row with zero
coefficients and a right-hand side excluding zero excludes the box.
Otherwise, back substitution bounds each multiplier.

Keep the nonnegative part of each multiplier interval and apply the
leaf upper bound where appropriate. Empty intervals exclude the box.
Using the retained intervals, evaluate
\[
 DF_G(\bc)^T\blambda-\nabla M(\bc),
 \qquad 2\sum_{e\in E}\lambda_e-M(\bc).
\]
Both expressions must vanish at a candidate minimum. Discard the box if any component excludes zero. The equality for the second moment follows
by pairing the vertex balances with the centred positions and using
the contact distances.

Every surviving box attempts root isolation after the prescribed
subdivision depth.
There are at most as many remaining contact equations as free angles.
When the remaining contact equations match the number of free angles,
they form the square system. When there are fewer equations, add
stationarity of the second moment in the angles, with one multiplier
per remaining contact. Tree contacts contribute nothing because their
distances are constant. The isolation box combines the angular box
with the computed bounds on these multipliers when needed.

The Krawczyk test computes an interval image of this box
\cite{krawczyk}. Disjointness from the box excludes roots.
Strict inclusion in its interior proves a unique root.
Otherwise subdivide. The isolation box contains the angles and
required multipliers of every possible minimum represented by the
angular box, including its boundary. The separate bridge-equation
test uses only the exclusion criterion.

Contract a unique-root enclosure and verify strict noncontacts before
retaining it as an \emph{isolated box} for candidate analysis.
The isolation equations impose the prescribed contacts but not strict
noncontacts, so a root may have additional contacts.
If a previously verified configuration with additional contacts lies
in the root box, uniqueness identifies it with that configuration.
Account for it under its larger contact graph
(Supplementary Notes, Section~6.3). Other inconclusive separation
tests cause subdivision.

Pending boxes are kept on a stack. Subdivision bisects the widest
angular interval and adds both closed children for the same tests.
Their union preserves coverage. Operational limits leave undecided
boxes as \emph{unresolved}. Search completion requires every box to
be excluded or accounted for by root isolation, with no unresolved boxes.
Comparisons of the second moment and second-order analysis follow in
candidate analysis.

\subsection{Candidate analysis}
\label{subsec:candidate-analysis}

The interval search completes all 12 and 44 assigned graphs for five
and six disks, respectively, with no unresolved boxes. The 9 and 33
graphs with leaves retain no candidate. The remaining searches retain
four isolated boxes for five disks and seven for six. By coverage and
the necessary condition tests, every possible minimum in a graph
assigned to the interval search is represented by a retained root,
up to isometry and relabelling
(Supplementary Notes, Proposition~10.1 and Section~10.2).

Each isolated box yields enclosures of the centred Cartesian coordinates
and contact multipliers. A contact minor of order \(|E|\) verifies the
rank hypothesis of Lemma~\ref{lem:constant-rank-contact-locus}.
With strict noncontacts already established, we evaluate the intrinsic Hessian using 
Proposition~\ref{prop:independent-defining-equations} and
\eqref{eq:contact-intrinsic-hessian}
\[
 \operatorname{Hess}_N M(\bc)[\bv,\bv]
 =\sum_{i=1}^n\|\bv_i-\overline{\bv}\|^2
 -\sum_{\{i,j\}\in E}\lambda_{ij}\|\bv_i-\bv_j\|^2.
\]
A negative value in a tangent
direction excludes local minimality by
Lemma~\ref{lem:second-order-exclusion}. A second moment above \(U_n\)
excludes global minimality.

All eleven retained boxes enclose roots with second moment strictly above the target.
Table~\ref{tab:retained-moment-bounds} gives bounds valid for every
retained box of each graph. Table~\ref{tab:candidate-verdicts} records
their exclusions, prioritising exclusion by a negative intrinsic
direction. Its numerical
columns count boxes, which need not represent distinct geometric
classes. Zero retained boxes means that the search leaves no candidate
with that contact graph; a zero exclusion count means that no
box was assigned to that reason.

\begin{table}[htbp]
\centering\small
\caption{Retained angular boxes, classified by exclusion reason.}
\label{tab:candidate-verdicts}
\begin{tabular}{@{}clrrr@{}}
\toprule
\(n\) & Graph6 & Retained boxes & Negative direction & Second moment bound\\
\midrule
5 & \texttt{D\textasciigrave{}\{} & 2 & 2 & 0\\
5 & \texttt{Dd[} & 1 & 1 & 0\\
5 & \texttt{DqK} & 1 & 1 & 0\\
\midrule
6 & \texttt{EDZW} & 1 & 1 & 0\\
6 & \texttt{ETXW} & 1 & 1 & 0\\
6 & \texttt{EoSo} & 1 & 1 & 0\\
6 & \texttt{EoSw} & 1 & 1 & 0\\
6 & \texttt{EElw} & 2 & 0 & 2\\
6 & \texttt{EQlw} & 1 & 0 & 1\\
6 & \texttt{E@ro} & 0 & 0 & 0\\
6 & \texttt{E@rw} & 0 & 0 & 0\\
6 & \texttt{EPVW} & 0 & 0 & 0\\
6 & \texttt{E\textasciigrave{}NW} & 0 & 0 & 0\\
6 & \texttt{E\textasciigrave{}ow} & 0 & 0 & 0\\
\bottomrule
\end{tabular}
\end{table}
\begin{table}[htbp]
\centering\small
\caption{Bounds valid for every retained box of each graph.
\(Q\) is the intrinsic quadratic value of a recorded tangent test
vector. For graphs with two boxes, the smaller lower bound for \(M\)
and, where available, the larger upper bound for \(Q\) are shown.
The decimals are rounded outwards from the rational interval records.
A dash means no internal tangent direction.}
\label{tab:retained-moment-bounds}
\begin{tabular}{@{}clrr@{}}
\toprule
\(n\) & Graph6 & Lower bound for \(M\) & Upper bound for \(Q\)\\
\midrule
5 & \texttt{D\textasciigrave{}\{} & 7.99 & \(-2.39\)\\
5 & \texttt{Dd[} & 6.98 & \(-0.34\)\\
5 & \texttt{DqK} & 7.23 & \(-4.87\)\\
\midrule
6 & \texttt{EDZW} & 10.23 & \(-0.60\)\\
6 & \texttt{EElw} & 9.99 & ---\\
6 & \texttt{EQlw} & 10.99 & ---\\
6 & \texttt{ETXW} & 11.46 & \(-8.61\)\\
6 & \texttt{EoSo} & 11.99 & \(-9.16\)\\
6 & \texttt{EoSw} & 11.50 & \(-6.13\)\\
\bottomrule
\end{tabular}
\end{table}

The retained roots for \texttt{EElw} and \texttt{EQlw} are rolling
rigid: their contact rank is nine, so no non-Euclidean tangent
direction is available for the second-order test.

The search for \texttt{EPVW} identified a root with an additional
contact, belonging to \texttt{EQlw}. An algebraic witness in the tested
box and root uniqueness account for this boundary case, so it leaves
no candidate with contact graph \texttt{EPVW}
(Supplementary Notes, Section~6.3).

The target graphs \texttt{DR\{} and \texttt{EDZw} are analysed
by computations in \(\Q(\sqrt3)\). Their graph6 vertex labels differ
from the coordinate labels in \eqref{eq:five-edges} and
\eqref{eq:six-edges}. If \(p(i)\) is the coordinate label corresponding
to graph6 vertex \(i\), then
\[
\begin{aligned}
n=5:\quad &(p(1),\ldots,p(5))=(1,3,4,5,2),\\
n=6:\quad &(p(1),\ldots,p(6))=(2,6,4,3,5,1).
\end{aligned}
\]
These relabellings identify the contact edges and their multipliers.
The target \texttt{EDZw} and the searched graph \texttt{EDZW} are
distinct graph classes. All two-circle attachment branches are
enumerated: \(4\) for five disks and \(8\) for six. Contact, strict
noncontact, congruence, and graph-automorphism checks leave one
realisation class for each target. Its \(M\), positive multipliers,
and full contact rank are verified (Supplementary Notes, Section~8).

The six disk graph \texttt{E\textasciigrave{}dg} has a
complete critical family on which \(M = 11\). At each member,
a tangent direction has negative intrinsic second variation, so the
family contains no ambient local minimum. Its completeness is proved
in Supplementary Notes, Section~9.

\subsection{Global classification}
\label{subsec:global-conclusion}

Let \(\Euc\) denote the full Euclidean isometry group of the plane,
and let \(S_n\) be the group of permutations of the disk labels.

\begin{theorem}[Global classification]
\label{thm:validated-classification}
For every \(\bc\in\D_5\), \(M(\bc)\geq34/5\), with equality only
on the \(\Euc\times S_5\) orbit of the two-row cluster
\eqref{eq:five-coordinates}. For every
\(\bc\in\D_6\), \(M(\bc)\geq29/3\), with equality only on the
\(\Euc\times S_6\) orbit of the missing-neighbour flower
\eqref{eq:six-coordinates}.
\end{theorem}
\begin{proof}
Choose a connected global minimiser by
Lemma~\ref{lem:connected-census-reduction}. By
Lemma~\ref{lem:outerplanar-census-reduction}, its contact graph belongs
to the census in Table~\ref{tab:proof-pipeline-summary}. The critical
family is excluded. Every possible minimum on a graph assigned to the
interval search lies in a retained root box. For each such root,
Table~\ref{tab:retained-moment-bounds} gives a lower bound for \(M\)
greater than \(U_n\). The minimiser
must therefore have the target graph. Branch enumeration leaves
one realisation class, with second moment \(U_n\).

It remains to exclude disconnected equality cases. The target is an
ambient strict local minimiser modulo Euclidean motions by
Corollaries~\ref{cor:five-local} and~\ref{cor:six-local}.
If a disconnected packing attained the global value, the last rotation
in the connectedness reduction would produce disconnected packings
of that same value converging to the connected target. They cannot be
congruent to it, contradicting strict local minimality.
\end{proof}

The computation uses \texttt{DiskSecondMomentProof}, available upon
request. The reported counts refer to the archived run
\texttt{20261029T063416Z-6281e07a4606}, performed at 512-bit precision
for two to six disks. The archive records the source revision, software
environment, reproduction parameters, graph inputs, graph conventions,
subdivision records and candidate analyses. Every search closes without unresolved branches, and every
retained root is excluded as a global minimiser.

\section*{Acknowledgements}
During the preparation of this article, we used OpenAI's
ChatGPT and Anthropic's Claude to explore mathematical arguments
guided by questions and conjectures formulated by the authors.
Through this interaction, preliminary mathematical arguments
and proof strategies were explored and refined, and potential
objections were examined. The assistance included the following
activities:

\begin{itemize}

\item Suggesting proof strategies and exploring alternative approaches
for the exact tangent cone formula
(Lemma~\ref{lem:exact-tangent-cone}), the centroidal contact
balance (Proposition~\ref{prop:ambient-first-order}), and
the local minimality criteria
(Theorem~\ref{thm:critical-cone} and
Corollary~\ref{cor:rigid-positive-balance}).

\item Exploring approaches to redundant contact constraints
and second order exclusion arguments.

\item Suggesting geometric reductions involving contact
leaves, exceptional configurations, and critical families.

\item Examining the global classification
(Theorem~\ref{thm:validated-classification}), particularly
the completeness of the contact graph census, boundary
cases, and the passage from local analysis to unrestricted
global optimality.

\item Assisting with the design and examination of
computational procedures, including contact graph enumeration,
root isolation, interval arithmetic, search coverage,
and candidate exclusions.

\item Identifying potential gaps in mathematical arguments,
suggesting corrections, and assisting with literature searches,
figures, LaTeX preparation, and editorial revisions.

\end{itemize}

The authors checked the mathematical arguments and computational
results incorporated into the final manuscript and take full
responsibility for its content.


\end{document}